\pdfoutput=1
\RequirePackage{ifpdf}
\ifpdf % We are running pdfTeX in pdf mode
\documentclass[pdftex]{sigma}
\else
\documentclass{sigma}
\fi

\newcommand{\rphi}{{\;_R\phi}}
\DeclareMathOperator{\id}{id}
\newcommand{\apn}{A//N}
\newcommand{\ccb}{\mathcal B}
\newcommand{\bdelta}{\Delta}
\newcommand{\bann}{\mtr{Ann}}
\newcommand{\ra}{\rightarrow}
\newcommand{\Z}{\mathbb Z}
\newcommand{\ot}{\otimes}
\newcommand{\co}{\mathcal O}
\newcommand{\mtc}{\mathcal}
\newcommand{\lam}{\lambda}
\newcommand{\Lam}{\Lambda}
\newcommand{\al}{\alpha}
\newcommand{\eps}{\epsilon}
\newcommand{\rh}{\rightharpoonup}
\newcommand{\lh}{\leftharpoonup}
\newcommand{\uw}{\uparrow}
\newcommand{\ch}{\chi}
\newcommand{\mtr}{\mathrm}
\newcommand{\ncm}{\newcommand}
\newcommand{\blam}{\Lam}
\ncm{\cc}{\mtc{C}}
\ncm{\Rep}{\mathrm{Rep}}
\ncm{\cA}{\mathcal{A}}
\ncm{\cd}{{\mathcal{D}}}
\ncm{\G}{\mathcal{G}}
\ncm{\cop}{\mtr{cop}}
\ncm{\lker}{\mtr{LKer}}
\ncm{\rep}{\Rep}
\newcommand{\irr}{\mathrm{Irr}}
\newcommand{\kk}{\Bbbk}
\newcommand{\lkera}{\lker_{A}}
\newcommand{\ca}{\cA}
\newcommand{\csu}{\mathbf{ C}}
\newcommand{\mtca}{\mtc A}
\newcommand{\apl}{A//L}
\newcommand{\dimka}{\dim_{\kk}(A)}
\newcommand{\dimk}{\dim_{\kk}}
\newcommand{\blaml}{\blam_L}
\newcommand{\dimkl}{\dim_{\kk}(L)}
\newcommand{\mtcjl}{\mtc J_L}
\newcommand{\ccad}{{\cc_{\mtr{ad}}}}
\newcommand{\ccpt}{{\cc_{\mtr{pt}}}}
\newcommand{\phir}{\phi_R}
\newcommand{\kda}{{ \Phi(A)}}
\newcommand{\mtcil}{\mtc{I}_L}

\numberwithin{equation}{section}

\newtheorem{Theorem}{Theorem}[section]
\newtheorem{Corollary}[Theorem]{Corollary}
\newtheorem{Lemma}[Theorem]{Lemma}
\newtheorem{Proposition}[Theorem]{Proposition}
{ \theoremstyle{definition}
\newtheorem{Example}[Theorem]{Example}
\newtheorem{Remark}[Theorem]{Remark} }

\begin{document}
\allowdisplaybreaks

\newcommand{\arXivNumber}{1709.02176}

\renewcommand{\PaperNumber}{039}

\FirstPageHeading

\ShortArticleName{Representations and Conjugacy Classes of Semisimple Quasitriangular Hopf Algebras}

\ArticleName{Representations and Conjugacy Classes \\ of Semisimple Quasitriangular Hopf Algebras}

\Author{Sebastian BURCIU}

\AuthorNameForHeading{S.~Burciu}

\Address{Institute of Mathematics ``Simion Stoilow'' of the Romanian Academy,\\ P.O.~Box 1-764, RO-014700, Bucharest, Romania}
\Email{\href{sebastian.burciu@imar.ro}{sebastian.burciu@imar.ro}}
\URLaddress{\url{http://www.imar.ro/~sburciu/}}

\ArticleDates{Received September 16, 2019, in final form April 27, 2020; Published online May 06, 2020}

\Abstract{In this paper we give two general formulae for the M\"uger centralizers in the category of representations of a semisimple quasitriangular Hopf algebra. The first formula is given in the terms of the Drinfeld map associated to the quasitriangular Hopf algebra. The second formula for the M\"uger centralizer is given in the terms of the conjugacy classes introduced by Cohen and Westreich in [\textit{J.~Algebra} \textbf{283} (2005), 42--62]. In the case of a factorizable Hopf algebra these formulae extend some particular cases obtained by the author in [\textit{Math.~Z.} \textbf{279} (2015), 227--240].}


\Classification{18D10; 16T05; 19D23}

\section{Introduction}

The notion of centralizer in a braided fusion category was introduced by M\"uger in \cite{proclond}. It was shown in \cite[Theorem~3.13]{dgno2} that the centralizer of a nondegenerate fusion subcategory of a~braided category is a categorical complement of the nondegenerate subcategory. This principle is the basis of many classification results of braided fusion categories, see for example papers \cite{DGNO, dgno2, eno-adv} and references therein.

Despite its importance, in the current literature there is no concrete known formula for the M\"uger centralizer of all fusion subcategories of a given fusion category. Only few particular cases are completely known in the literature. For instance, in the same paper~\cite{proclond}, M\"uger described the centralizer of all fusion subcategories of the category of finite-dimensional representations of a Drinfeld double of a finite abelian group. More generally, for the category of representations of a (twisted) Drinfeld double of an arbitrary finite group, not necessarily abelian, a similar formula was then given in~\cite{nnw}. For the braided center of Tambara--Yamagami categories, in~\cite{gnn}, the centralizer can be described by computing completely the $S$-matrix of the modular category. In~\cite{bcg} a different approach gave a partial formula for the centralizer of fusion subcategories of a~braided equivariantized fusion category.

Given a fusion subcategory $\cd$ of a braided fusion category~$\cc$, the notion of {\it M\"uger centralizer of~$\cd$} was introduced in~\cite{dgno2}. The centralizer~$\cd'$ is defined as the fusion subcategory $\cd'$ of~$\cc$ generated by all simple objects $X$ of $\cc$ satisfying
\begin{gather*} c_{X, Y}c_{Y, X}=\mtr{id}_{X\ot Y}\end{gather*}
for all objects $Y\in\co(\cd)$ (see also~\cite{proclond}). For a fusion category~$\cc$ as usually, we denote by $\co(\cc)$ the set of isomorphism classes of simple objects of~$\cc$.

If $(A,R)$ is a quasitriangular Hopf algebra then the category $\rep(A)$ of finite-dimensional $A$-modules is a braided category with the braiding given by
\begin{gather*}
c_{M, N}\colon \ M\ot N\ra N\ot M,\qquad m\ot n \mapsto R_{21}(n\ot m)=R^{(2)}n\ot R^{(1)}m,
\end{gather*}
for any two objects $M, N\in \rep(A)$.

Given a quasitriangular Hopf algebra $(A, R)$ one can also define the Drinfeld map
\begin{gather*}
\phi_R\colon \ A^* \ra A, \qquad f \mapsto (f \ot \id)(R_{21}R)=f(Q_{1})Q_{2},
\end{gather*}
where $Q=R_{21}R$ is the monodromy matrix.

We prove the following theorem which gives a general description for the centralizer of any fusion subcategory of the category of representations of a quasitriangular Hopf algebra:
 \begin{Theorem} \label{main1}Let $(A, R)$ be a semisimple quasitriangular Hopf algebra and $L$ be a left normal coideal subalgebras of $A$. Then \begin{gather*} \rep(A//{L})'= \rep(A//{M}),\qquad \text{where}\quad {M}=\phi_{R}((A//{L})^{*}).
\end{gather*}
\end{Theorem}
We denote by $F_{0}, F_{1}, \dots, F_{r}$ the central primitive idempotents of the character ring~$C(A)$ where $F_{0}=t$ is the idempotent integral of~$A^{*}$. Following~\cite{CW2} one can define the conjugacy classes~$\cc^{j}$ of~$A$ as $\cc^{j}:=\blam \lh F_{j}A^{*}$, where $\blam$ is an idempotent integral of $A$ and $a \lh f=\langle f, a_{1}\rangle a_{2}$ for all $a \in A$ and $f \in A^{*}$. It is well known that these conjugacy classes are the simple $D(A)$-submodules of the induced $D(A)$-module $\kk\uparrow^{D(A)}_{A}\simeq A$, see~\cite{zind}.

Let $(A, R)$ be a semisimple quasitriangular Hopf algebra and $V_{0}=\kk, \dots ,V_{r}$ be a complete set of isomorphism classes of irreducible $A$-modules. Let also $\irr(A)=\{\ch_{0}=\eps,\, \ch_{1}, \dots ,\ch_{r}\}$ be the set of irreducible characters afforded by these modules and $E_{i}\in \mtc Z(A)$ be the associated central primitive idempotent of the irreducible character~$\ch_{i}$. Since the Drinfeld map $\phi_{R}\colon C(A)\ra \mtc Z(A) $ is an algebra map we may suppose that $\phi_{R}(F_{j})=\sum\limits_{i \in \mtc A_{j}}E_{i}$ for some subset $\ca_j\subseteq \{0, \dots , r\}$. Since $\phi_{R}(1)=1$ we obtain a partition for the set of indices of all irreducible representations $\{0,1, \dots, r\}= \bigsqcup\limits_{j\in \mtc J}\mtc A_{j}$. For any $0\leq i\leq r$ we denoted by $m(i)$ the unique index $j \in \mtc J$ such that $i \in \mtc A_{j}$. Therefore in this way we obtain a unique function
\begin{gather*}
m\colon \ \{0, 1,\dots, r\}\ra \mtc J
\end{gather*}
with the property that $E_i\phir(F_{m(i)})\neq 0$ for all $i \in \{0, 1,\dots, r\}$.

Our second main result is the following:
\begin{Theorem} \label{mgqtr}Suppose that $(A,R)$ is a semisimple quasitriangular Hopf algebra and $L$ is a left normal coideal subalgebra of~$A$. With the above notations one has
\begin{gather*}
\co(\rep(A//L)')=\big\{\ch_{i}\,|\, \cc^{m(i)}\subseteq L\big\}.
\end{gather*}
\end{Theorem}
Recall that the quasitriangular Hopf algebra $(A,R)$ is called {\it factorizable} if the Drinfeld map $\phi_{R}\colon A^*\ra A$ is an isomorphism of algebras. In this case, its restriction $\phir|_{C(A)}\colon C(A)\ra \mtc Z(A)$ is an isomorphism of algebras. For a factorizable semisimple Hopf algebra we can record the primitive central idempotents $F_j$ of $C(A)$ such that $F_{j}:=\phi_{R}^{-1}(E_{j})$ any $1\leq j \leq r$. With these notations, $m(i)=i$ for any $0\leq i\leq r$ and Theorem~\ref{mgqtr} implies the following:
\begin{Corollary}\label{main2}
Let $(A, R)$ be a semisimple factorizable Hopf algebra and~$L$ be a left normal coideal subalgebra of~$A$. Then with the above notations one has that
\begin{gather*}
\co(\rep(A//L)')=\big\{\ch_{i}\,|\,\cc^{i}\subseteq L\big\}.
\end{gather*}
\end{Corollary}

Shortly, this paper is organized as follows. In Section~\ref{prelim} we recall the basic notions of Hopf algebras and fusion categories that are used throughout this paper. In this section we also prove a canonical decomposition of a left normal coideal subalgebra in terms of the decomposition of its integral, see equation~\eqref{decl3}. In Section~\ref{qtr} we recall the main properties of quasitriangular Hopf algebras and their associated Drinfeld maps. In Section~\ref{fmr} we prove Theorem~\ref{main1} and some consequences of it. In particular we apply Theorem~\ref{main1} to the adjoint subcategory of the category of representations of a factorizable Hopf algebra. In this way we obtain a relation, via the Drinfeld map, between the Hopf center and the first commutator of a factorizable Hopf algebra.

In Section~\ref{smr} we prove Theorem~\ref{mgqtr}. Some consequences of this result are also described. In Section~\ref{h8} we give an example, by considering the semisimple quasitriangular Hopf algebra~$H_8$ of dimension~$8$. Based on our results we are able to compute the function~$m$ in this case and therefore the centralizer of any fusion subcategory of~$\rep(H_8)$.

We work over an algebraically closed field~$\kk$ of characteristic zero. The comultiplication and antipode of a Hopf algebra are denoted by~$\Delta$ and~$S$ respectively. We use Sweedler’s notation for comultiplication with the sigma symbol dropped. All the other Hopf algebra notations are those used in~\cite{montg}.

\section{Preliminaries}\label{prelim}
Let $A$ be a finite-dimensional semisimple Hopf algebra over an algebraically closed field $\kk$ of characteristic zero. Then~$A$ is also cosemisimple and $S^{2}=\id$~\cite{Lard}. The character ring
$C(A):=\G_0(A)\ot_{\Z}\kk$ is a semisimple subalgebra of~$A^*$ and it has a vector space basis
given by the set $\mtr{Irr}(A)$ of irreducible characters of~$A$, see~\cite{Z}. Moreover,
$C(A)=\mtr{Cocom}(A^*)$, the space of cocommutative elements of~$A^*$. By duality, the
character ring of $A^*$ is a semisimple subalgebra of $A$ and $C(A^*)=\mtr{Cocom}(A)$. If
$M$ is an $A$-representation with character $\chi$ then~$M^*$ is also an
$A$-representation with character $\chi^*=\chi \circ S$. This induces an involution
``$\;^*\;$'': $C(A)\ra C(A)$ on~$C(A)$.

Throughout of this paper we denote by $\blam$ an idempotent integral of $A$ and by $t$ an idempotent integral of $A^{*}$. Moreover one has that $t(\blam)=\frac{1}{\dim_{\kk}(A)}$. Recall also~\cite{Lar} that
\begin{gather*}
\dim_{\kk}(A)\blam=\sum_{d \in \irr(A^{*})}\eps(d)d
\end{gather*}
is the regular character of $A^{*}$. Dually
\begin{gather*}
\dim_{\kk}(A)t=\sum_{\ch \in \irr(A)}\ch(1)\ch
\end{gather*}
is the regular character of~$A$.

Recall that a {\it left} coideal subalgebra of~$A$ is a subalgebra $L\subseteq A$ with $\bdelta(L)\subset A\ot L$. Then~$L$ is called {\it left normal coideal subalgebra} if $L$ is closed under the left adjoint action of $A$, i.e., $a_{1}lS(a_{2})\in L$ for any $l \in L$ and any $a\in A$. Recall also from~\cite{gmj} that given a left coideal subalgebra~$L$ of~$A$ there is a unique element $\blam_{L}\in L$ (called {\it integral}) such that $l\blam_{L}=\eps(l)\blam_{L}$ for all $l \in L$, see also~\cite{kopp}. Then the coideal subalgebra~$L$ is normal if and only if~$\blam_{L}$ is central, i.e., $\blam_L\in Z(A)$.

For any left normal coideal subalgebra $L$ of $A$ the augmentation ideal $AL^{+}$ is a Hopf ideal and it has the following form $AL^{+}=A(1-\Lam_{L})=\bann_{A}(\Lam_{L})$. Thus one can define the Hopf quotient $A//L:=A/AL^+$. It is well known that any fusion subcategories of~$\rep(A)$ can be written as $\rep(A//L)$ for some Hopf quotient of~$A$.
\begin{Remark}\label{regq} Since $A$ is free as left $L$-module \cite{sk} it follows that the map
\begin{gather*}%\label{regcid}
A\ot_{L}\kk\simeq A\blam_{L},\qquad a\ot_{L}1\mapsto a\blam_{L}
\end{gather*}
is an isomorphism of $A$-modules. Moreover, by \cite[Proposition~3.11]{iop} it follows that the regular module of the quotient Hopf algebra~$A//L$ is isomorphic to the induced module~$A\ot_{L}\kk$.
\end{Remark}
\subsection{Duality between the character ring and the center}
Let $A$ be a semisimple Hopf algebra over the ground field $\kk$. Let us denote by $\irr(A)$ the set of irreducible characters of $A$. We suppose that $\irr(A)=\{\ch_{0}, \ch_{1}, \dots, \ch_{r}\}$. Without loss of generality we may suppose that $\ch_{0}=\eps$. Let also $E_{0}, E_{1}, \dots, E_{r}$ be the corresponding central primitive central idempotents in $A$.
The evaluation form
\begin{gather*}\label{evform}
C(A)\otimes \mtc Z(A)\ra \kk,\qquad \ch \otimes a\mapsto \ch(a)
\end{gather*}
is nondegenerate. A pair of dual bases for this form is given by $\big\{\ch_{i}, \frac{1}{n_{i}}E_{i}\big\}$
since $\big\langle \ch_{i}, \frac{1}{n_{j}}E_{j}\big\rangle =\delta_{i,j}$ for any $1 \leq i,j \leq r$.

\subsection{Another pair of dual basis in the commutative case} Let $A$ be a semisimple Hopf algebra with a commutative character ring. According to~\cite{CW2} in the case of a commutative ring~$C(A)$ there is another pair of dual bases corresponding to this nondegenerate form. This pair of dual bases is given in terms of the conjugacy class sums as defined in~\cite{CW2}.

The {\it conjugacy classes $\cc^{j}$} of $A$ are defined as $\cc^{j}=\blam \lh F_{j}A^{*}$, where $\blam=\blam_A$ is a two-sided idempotent integral of $A$ and $\{F_j\}_j$ is the (complete) set of central primitive idempotents of the semisimple algebra $C(A)$. This notion of conjugacy classes generalizes the usual notion of conjugacy classes in finite groups.

\begin{Example}Let $G$ be a finite group and $A=\kk G$ be the associated group algebra. It is easy to see that the conjugacy classes as defined above coincide with the usual notion of conjugacy class in a group. Indeed, let ${\mtc C}_0, {\mtc C}_1, \dots ,\mtc C_r$ be the usual conjugacy classes of~$G$.
Then the set of central primitive idempotents of $C(\kk G)$ can be described as $p_{j}=\sum\limits_{h \in \cc^j}p_h$ where $p_h\in \kk^G$ is defined as $p_h(g)=\delta_{g,h}$. Since $\blam=\frac{1}{|G|}\sum\limits_{g \in G}g$ it follows that
\begin{gather*}
\cc^j=\blam\lh p_j\kk^G=\kk [\mtc C_j]
\end{gather*}
is the vector sub-space of $\kk G$ generated by all group elements of~$\cc_j$.
\end{Example}

Recall that the Fourier transform $\mtc{F}\colon A^*\ra A$ defined by $f\mapsto \blam \lh f$ for any $f \in A^*$ is a~$\kk$-linear isomorphism. Since $A^*=\bigoplus\limits_{j=0}^rF_jA^*$ and $\cc^j=\mtc{F}(F_jA^*)$ one has
\begin{gather*}
A=\bigoplus_{j=0}^{r}{\cc^{j}}.
\end{gather*}

One can also define the corresponding {\it conjugacy class sums}
\begin{gather*}%\label{clssum}
{\bf C}_{j}=\blam \lh (\dim A)F_{j}.
\end{gather*}
Note that ${\bf C}_j\in \mtc Z(A)$ and since $\dim_{\kk}\mtc Z(A)=\dim_{\kk}C(A)$ it follows that ${\cc^{j}}\cap \mtc Z(A)=\kk {\bf C}_j$. Note also that $\dimka\blam=\sum\limits_{j=0}^{r}{\bf C}_{j}$. Indeed, $\blam=\blam \lh \eps=\blam\lh \sum\limits_{j=0}^rF_j= \frac{1}{\dimka}\sum\limits_{j=0}^r\csu_j$.
\begin{Example}

If $A=\kk G$ then
\begin{gather*}
{\bf C}_{j}=\bigg(\frac{1}{|G|}\sum_{g \in G}g\bigg) \lh |G|p_{j}=\sum_{g \in G}p_j(g)g=\sum_{h \in C_j}h
\end{gather*}
is the usual class sum of a conjugacy class $\mtc C_j$.
\end{Example}
\begin{Remark}\label{genmz} By the class equation for semisimple Hopf algebras, see~\cite{leq}, one has that the value $n_{j}:=\frac{\dim_{\kk}A^{*}}{\dim_{\kk}(A^{*}F_{j})}$ is an integer. Moreover as in \cite[equation~(11)]{CW2} one can deduce that $F_{j}(\Lam)=\frac{1}{n_{j}}$.
\end{Remark}
\begin{Example}If $A=\kk G$ then $n_j=\frac{|G|}{|C_j|}$ is the order of the centralizer of any group element $g_j \in \mtc C_j$.
\end{Example}

This implies that a second pair of dual bases for the form of equation~\eqref{evform} can be given by $
\big\{F_{i}, \frac{n_{i}}{\dim_{\kk}(A)}{\bf C}_{i}\big\}$, see also \cite[equation~(17)]{CW2}. This can be written as $\big\langle F_{i}, \frac{n_{j}}{\dim_{\kk}(A)}{\bf C}_j\big\rangle =\delta_{i,j}$.

Indeed, one has
\begin{gather*}
\left\langle F_{i}, \frac{n_{j}}{\dim_{\kk}(A)}{\bf C}_j\right\rangle = \frac{n_{j}}{\dim_{\kk}(A)}\langle F_i, \blam\lh (\dimka)F_j\rangle
 = {n_{j}}\langle F_jF_i, \blam\rangle =\delta_{i,j}n_jF_j(\blam)=\delta_{i,j}.
\end{gather*}

\subsection{Decomposition of the integral}
Let $L$ be a left normal coideal subalgebra of a semisimple Hopf algebra $A$ with a commutative character ring~$C(A)$. We shall use the notation $\lam_{L}\in (A//L)^{*}$ for the idempotent integral of the Hopf algebra $(A//L)^{*}$. Clearly $\lam_{L}\in C((A//L)^{*})\subset C(A^{*})$ and we may suppose that
\begin{gather}\label{intldec}
\lam_L=\sum_{j \in \mtcil}F_j
\end{gather}
for some subset of indices $\mtcil\subseteq \{0,1, \dots, r\}$. Note that by \cite[Lemma~1.1]{CW10}
$\widetilde{\blam_{L}}:=\blam \lh \lam_{L}$ is a~left integral for~$L$. It follows from above that
\begin{gather*}
\widetilde{\blam_{L}}=\frac{1}{\dimka}\sum_{j \in \mtcil}\csu_j.
\end{gather*}
Then one has \begin{gather*} %\label{decompl}
L=\blam_{L}\lh A^{*}=\bigoplus_{j \in \mtc I_{L}}{\bf C}_j\lh A^{*}=\bigoplus_{j \in \mtc I_{L}}{\cc^{j}}. \end{gather*}
Note also that \begin{gather*}
\eps(\csu_j)=
\dimka F_j(\blam)=\frac{\dimka}{n_j}=\dimk(F_jA^*)=\dimk\big(\cc^j\big).
\end{gather*}
Thus
\begin{gather*}
\eps\big(\widetilde{\blam_{L}}\big)=\frac{1}{\dimka}\bigg(\sum_{j \in \mtcjl}\eps(\csu_j)\bigg)=\frac{1}{\dimka}\bigg(\sum_{j \in \mtcjl}\dimk\big(\cc^j\big)\bigg)=\frac{\dimkl}{\dimka}.
\end{gather*}
It follows then that{\samepage
\begin{gather}\label{decl3}
\blaml=\frac{\dimka}{\dimkl}\widetilde{\blam_{L}}=\frac{1}{\dimkl}\sum_{j \in \mtcil}\csu_j
\end{gather}
 is a formula for the idempotent integral of~$L$.}

Define the functional $p_{{\cc^{j}}}\in A^{*} $ as the unique linear functional that coincides to~$\eps$ on~${\cc^{j}}$ and it is equal to zero on the other conjugacy classes $\mtc C_{l}$ with $l \neq j$. The following lemma was proven in \cite[Theorem~5.13]{repalg}.

\begin{Lemma} \label{fct}Suppose that $A$ is a semisimple Hopf algebra with a commutative character ring $C(A)$. Let $\{F_j\}_{0\leq j\leq r}$ be a complete set of central primitive idempotents of~$C(A)$. Then $F_{j}=p_{{\cc^{j}}}$ for all $0\leq j\leq r$.
\end{Lemma}

Equation \eqref{decl3} and the above lemma implies the following:
\begin{Lemma}\label{charval} Let $L$ be a left normal coideal subalgebra of a semisimple Hopf algebra~$A$. With the above notations one has $j \in \mtc I_{L} \iff F_{j}(\blam_{L})\neq 0$.
\end{Lemma}

Let $A$ be a semisimple Hopf algebra with a commutative character ring $C(A)$. Then $\{F_{j}\}$ form a $\kk$-linear basis for $C(A)$ and for any character $\ch \in C(A)$ one can write $\ch=\sum\limits_{j=0}^{r}\al_{\ch,j}F_{j}$ with $\al_{\ch,j }\in \kk$.
Previous lemma implies the following result, see also \cite[Theorem~1.12]{CW-r}.

\begin{Proposition} \label{charl}Let $A$ be a semisimple Hopf algebra with a commutative character ring and $\ch \in C(A)$. Then one has
$ \ch \in C(A//L) \iff \ch F_{j}=\ch(1)F_{j}$ for all $j \in \mtc I_{L}.$
\end{Proposition}

\begin{proof} We may suppose that $\ch=\ch_{M}$ is the character of an $A$-module $M$. If $\ch \in C(A//L)$ then $\ch\lam_{L}=\ch(1)\lam_{L}$ and equation~\eqref{intldec} implies that $\ch F_{j}=\ch(1)F_{j}$, for any $j \in \mtc J_{L}$. Conversely if $\ch F_{j}=\ch(1)F_{j}$ for all $j \in J_{L}$ then $\ch \lam_{L}=\lam_{L}\ch(1)$. Thus
\begin{gather*}
\ch(\blam_{L})=\frac{\dimka}{\dimkl}\ch\big(\widetilde{\blam_{L}}\big)=\frac{\dimka}{\dimkl}\ch(\lam_{L}\rh \blam)=\frac{\dimka}{\dimkl}\ch\lam_{L}(\blam)=\frac{\dimka}{\dimkl}\ch(1)\lam_{L}(\blam).\!
\end{gather*}
On the other hand note that $\lam_{L}(\blam)=\eps\big(\widetilde{\blaml}\big)=\frac{\dimkl}{\dimka}
$ by equation~\eqref{decl3}. Thus $\ch(\blam_{L})=\ch(1)$ which shows that the restriction of the $A$-module $M$ to~$L$ is trivial. It follows that $M\in \rep(A//L)$.
\end{proof}
\begin{Example}If $A=\kk G$ then a left normal coideal subalgebra of $A$ is of the form $L=\kk N$ for some normal subgroup $N\unlhd G$. Then
\begin{gather*}
\lam_L=\sum_{C_j\subseteq N}p_j \qquad \text{and}\qquad \mtc{I}_L=\{j\,|\, C_j\subseteq N\}.
\end{gather*}
\end{Example}

\subsection{Left kernels and Burnside formula}\label{lkern}
Let $M$ be an $A$-module and let $\mtr{LKer}_{ _A}(M)$ be the {\it left kernel of $M$}. Recall \cite{gmj} that $\mtr{LKer}_{ _A}(M)$ is defined by
\begin{gather*}%\label{L}
\mtr{LKer}_{ _{A}}(M)=\{a \in A\,|\, a_1\ot a_2m=a\ot m,\,\text{for all}\, m\in M\}.
\end{gather*}
 Then by~\cite{gmj} it follows that $\lker_A(M)$ is the largest left coideal subalgebra of $A$ that acts trivially on~$M$. It is also a left normal coideal subalgebra. For example, if $A=\kk G$ is the group algebra of a finite group $G$ and~$M$ is a $\kk G$-module then $\lker_A(M)=\kk \ker_G(M)$, where $\ker_G(M)=\{g \in G\,|\, g.m=m, \,\text{for all}\, m\in M\}$ is the usual kernel of~$M$.

Next theorem generalizes a well known result of Brauer in the representation theory of finite groups.
\begin{Theorem}[{\cite[Theorem~4.2.1]{gmj}}]\label{charofim} Suppose that $M$ is a finite-dimensional module over a~semisimple Hopf algebra~$A$. Then \begin{gather*}
\langle M\rangle =\mtr{Rep}(A//\mtr{LKer}_A(M)),
\end{gather*}
where $\langle M\rangle $ is the fusion subcategory of~$\mtr{Rep}(A)$ generated by $M$.
\end{Theorem}

This implies that for any left normal coideal subalgebra $L$ of $A$ one has that
\begin{gather}\label{brint}
\bigcap_{M\in \irr(A//L)} \lker_{A}(M)=L.
\end{gather}

The previous theorem also implies that any fusion subcategory of $\rep(A)$ is of the type $\rep(A//L)$ for some left normal coideal subalgebra~$L$ of~$A$. Moreover, for any $V\in \rep(A)$ one has
 \begin{gather}\label{inq}
M\in \rep(A//L)\iff \lkera(M)\supseteq L.
\end{gather}

\subsection{Quasitriangular and factorizable Hopf algebras}\label{qtr}

Recall that a Hopf algebra $A$ is called {\it quasitriangular} if $A$ admits an $R$-matrix, i.e., an element $R \in A\ot A$ satisfying the following properties:
\begin{enumerate}\itemsep=0pt
\item[1)] $R \Delta(x)=\Delta^{\cop}(x)R $ for all $x \in A$,
\item[2)] $(\Delta \otimes \id)(R)= R^{(1)} \ot r^{(1)} \ot R^{(2)}r^{(2)}$,
\item[3)] $( \id \otimes \Delta)(R)=R^{(1)}r^{(1)} \ot r^{(2)}\ot R^{(2)}$,
\item[4)] $( \id \otimes \eps)(R)=1=(\eps \ot \id)(R)$.
\end{enumerate}
Here $R=r=R^{(1)}\ot R^{(2)}=r^{(1)}\ot r^{(2)}$. If $(A, R)$ is a quasitriangular Hopf algebra then the category of representations is a braided fusion category with the braiding given by
\begin{gather*}
c_{M, N}\colon \ M\ot N\ra N\ot M,\qquad m\ot n \mapsto R_{21}(n\ot m)=R^{(2)}n\ot R^{(1)}m
\end{gather*}
for any two left $A$-modules $M, N\in \rep(A)$ (see~\cite{Kas}). Recall that $R_{21}:= R^{(2)}\ot R^{(1)}$.
Denote $Q:=R_{21}R$. Then the monodromy of two objects is defined as
\begin{gather*}%\label{mdrmy}
c_{M,N}c_{N,M}\colon \ M\ot N \ra N\ot M,\qquad (m\ot n)\mapsto R^{(2)}R^{(1)}m\ot R^{(2)}R^{(1)}n=Q(m\ot n).
\end{gather*}

A quasitriangular Hopf algebra $(A, R)$ is called {\it factorizable} if and only if the Drinfeld map
\begin{gather*}%\label{dr}
\phi_R\colon \ A^* \ra A, \qquad f \mapsto (f \ot \id)(R_{21}R)
\end{gather*}
is an isomorphism of vector spaces. In this situation, following \cite[Theorem~2.3]{schfact} $\phi_R$ maps the character ring~$C(A)$ onto the center $\mtc Z(A)$ of~$A$ and the restriction $\phir|_{C(A)}$ is an isomorphism of algebras.
\begin{Remark} \label{nat} By \cite[Lemma~4.1]{rqts} one has that $\phi_{R}(C)$ is a left normal coideal for any subcoalgebra $C$ of $A^{*}$.
\end{Remark}

One can also define the map $\rphi\colon A^*\ra A$ by $ _{R}\phi(f)=(\id \ot f)(Q)$ for all $f\in A^*$. Moreover by \cite[Theorem~2.1]{schfact} one has that ${}_R\phi(f\chi)={}_R\phi(f)\,{}_R\phi(\chi)$ for all $f \in A^*$ and $\chi \in C(A)$. Thus $_{R}\phi|_{C(A)}\colon C(A)\ra \mtc Z(A)$ is an isomorphism of $\kk$-algebras.

\begin{Remark}\label{rphir} By \cite[Lemma~2.3]{rqts} one has that $S\circ_{R}\phi=\phi_{R}\circ s$ where $S$ and $s$ are the antipodes of $A$ and $A^{*}$ respectively.
\end{Remark}

In the case of a factorizable Hopf algebra $_{R}\phi$ is also bijective and moreover by~\cite{CW5} the two maps $_{R}\phi$ and $\phir$ coincide on the character ring~$C(A)$.

\section{Proof of Theorem \ref{main1} on M\"uger centralizer}\label{fmr}
 In this section we prove the first main theorem mentioned in the introduction. Given a fusion subcategory $\cd$ of a braided fusion category~$\cc$, recall that the M\"uger centralizer~$\cd'$ is defined as the fusion subcategory of~$\cc$ generated by all simple objects~$X$ of~$\cc$ satisfying
\begin{gather*} c_{X, Y}c_{Y, X}=\mtr{id}_{X\ot Y}\end{gather*}
for all objects $Y \in \co(\cd)$ (see also~\cite{proclond}). Recall that $\co(\cd)$ denotes the set of isomorphism classes of simple objects of~$\cd$.

Let $A$ be a semisimple quasitriangular Hopf algebra over $\kk$ and $\cd=\rep(A//L)$ be a fusion subcategory of $\rep(A)$ where $L$ is a left normal coideal subalgebra of~$A$.

\begin{Lemma}\label{firstequiv}
Let $(A, R)$ be a semisimple quasi\-trian\-gular Hopf algebra and $L$, $N$ be two left normal coideal subalgebras. Then the following assertions are equivalent:
\begin{enumerate}\itemsep=0pt
\item[$1)$] $\rep(A//N)\subseteq \rep(A//L)'$,
\item[$2)$] the following equation holds in $A\ot A$:
\begin{gather}\label{cond1}
Q(\blam_{L}\ot \blam_{N})=\blam_{L}\ot \blam_{N},
\end{gather}
\item[$3)$] $N\supseteq \phi_R((A//L)^*)$.
\end{enumerate}
 \end{Lemma}

\begin{proof}$(1)\iff (2)$ It is well known that two fusion subcategories of $\rep(A)$ centralize each other if and only if their regular representations centralize. Thus one needs to show that the two regular representations of $A//L$ and $A//N$ centralize each other if and only if equation~\eqref{cond1} holds. On the other hand from the definition of the braiding in $\rep(A)$ the two representations centralize each other if and only if $Q=R^{(2)}r^{(1)}\ot R^{(1)}r^{(2)}$ acts as identity on their tensor product $A//L\ot A//N$. By Remark~\ref{regq} one has $\kk\uw^{A}_{L}\ot \kk\uw^{A}_{N}=A\blam_{L}\ot A\blam_{N}$. Since $\blam_{L}$ and $\blam_{N}$ are central elements of $A$ it is clear that~$Q$ acts as identity on this subspace of $A\ot A$ if and only if equation~\eqref{cond1} holds.

$(2)\implies (3)$ By \cite[Lemma~1.1]{CW10} one has $(\apl)^*=\blam_L\rh A^*$. Therefore $\phi_{R}(\blam_{{L}}\rh f)=f(Q^{1}\blam_{L})Q^{2}$ for any $f \in A^{*}$. From here, applying equation~\eqref{cond1} it follows
\begin{gather}\label{condx}
\phi_{R}(\blam_{{L}}\rh f)\blam_{{N}}=f\big(Q^{1}\blam_{{L}}\big)Q^{2}\blam_{{N}}=f(\blam_{{L}})\blam_{{N}}
.\end{gather}
On the other hand note that \begin{gather}\label{inp} \eps_A(\phi_{R}(\blam_{{L}}\rh f))=f(\blam_{{L}}). \end{gather}
If $L':=\phi_{R}((A//{L})^{*})$ then equation~\eqref{inp} gives that $A(L')^{+} \subseteq AN^{+}$ and by \cite[Lemma~6.2]{iop} one has that $L'\subseteq {N}$.

$(3)\implies (2)$ In this situation equation~\eqref{condx} is satisfied for any $f \in A^*$ which shows that equation~\eqref{cond1} also holds.\end{proof}
\begin{Remark}
If $\ccb\subseteq \cd$ are fusion subcategories of a braided fusion category $\cc$ then clearly $\cd'\subseteq \ccb'$.
In particular if $
 \rep(A//N)\subseteq \rep(A//L)'$ then by centralizing once more one has that $\rep(A//L)\subseteq \rep(A//L)''\subseteq \rep(\apn)'$. Thus the three above conditions from the previous lemma are also equivalent to:
\begin{enumerate}\itemsep=0pt
\item[1)] $\rep(A//L)\subseteq \rep(\apn)'$,
\item[2)] $Q(\blam_{N}\ot \blam_{L})=(\blam_{N}\ot \blam_{L})$,
\item[3)] $L\supseteq \phi_R((A//N)^*)$.
\end{enumerate}
\end{Remark}

\subsection{Proof of Theorem \ref{main1}}

\begin{proof} Let $M:=\phir((A//L)^*)$ as in the statement of the theorem. By Remark~\ref{nat} it is well-known that~$M$ is also a left normal coideal subalgebra of $A$. Suppose also that $\rep(\apl)'=\rep(A//L^\circ)$ for some other left normal coideal subalgebra $L^\circ$ of $A$.

We need to prove that $L^\circ=M$. Note that for any left normal coideal subalgebra $N$ of $A$ one has $\rep(\apn)\subseteq \rep(\apl)'$ if and only if $\rep(\apn)\subseteq \rep(A//L^\circ)$, i.e $L^\circ\subseteq N$. Then the previous Lemma~\ref{firstequiv} shows that for any left normal coideal subalgebra~$N$ of~$A$ one has
\begin{gather*}
L^\circ\subseteq N\iff M\subseteq N.
\end{gather*}
In particular, for $N=L^\circ$ one obtains that $M\subseteq L^\circ$. For $N=M$ one obtains the other inclusion $L^\circ\subseteq M$. Thus $L^\circ=M$ and the proof is complete.
\end{proof}

Theorem~\ref{main1} can be rewritten as following by using the notion of left kernel of an $A$-module, see Section~\ref{lkern}:
\begin{Theorem}\label{putgh1} Let $(A, R)$ be a quasitriangular semisimple Hopf algebra and~$L$ be a left normal coideal subalgebra of~$A$. If $M\in \rep(A)$ then the following assertions are equivalent:
\begin{enumerate}\itemsep=0pt
\item[$1)$] $M\in \rep(A//L)'$,
\item[$2)$] $\phi_{R}((A//L)^{*})\subseteq \lker_{A}(M)$.
\end{enumerate}
\end{Theorem}
\begin{proof}One has that $\rep(A//L)'=\rep(A//L^\circ)$ where $L^\circ=\phi_{R}((A//L)^{*})$. Then by equation~\eqref{inq} one has $M\in \rep(A//L^\circ)\iff \lkera(M)\supseteq L^\circ$.
\end{proof}

\subsection{On the commutators and Hopf centre}
Given a fusion category~$\cc$ we denote by $\ccpt$ the fusion subcategory generated by the invertible objects of~$\cc$. In the case $\cc=\rep(A)$ for a semisimple Hopf algebra we have that $\ccpt$ is the full abelian subcategory generated by one-dimensional modules. It was shown in~\cite{iop} that $\rep(A)_{\rm pt}=\rep(A/I)$ where $I:=\{ab-ba\,|\, a, b\in A\}$ is the first commutator of the $\kk$-algebra~$A$. Moreover the commutator ideal $[A, A]$ is a Hopf ideal and by Takeuchi's correspondence it corresponds to a left normal coideal subalgebra~$A'$. Thus $A(A')^+=I$ and
\begin{gather*}
\rep(A)_{\rm pt}=\rep(A//A').
\end{gather*}
Moreover by \cite{iop} $A'$ is the smallest left normal coideal subalgebra~$L$ with the property that~$A//L$ is a commutative Hopf algebra. $A'$ is called the {\it the commutator} of~$A$.
\begin{Example}
If $A=\kk G$ then $A'=\kk G'$ where $G'$ is the first commutator of $G$, i.e., $G'=[G,\;G]$.
\end{Example}

For a fusion category $\cc$, recall that the adjoint subcategory $\cc_{\rm ad}$ is defined as the smallest fusion subcategory generated by all objects of the type $X\ot X^*$ with $X$ a~simple object of~$\cc$. If $\cc=\rep(A)$ is the category of representations of a semisimple Hopf algebra~$A$ then it is well-known that
\begin{gather*}
\rep(A)_{\rm ad}=\rep(A//K(A)),
\end{gather*}
where $K(A)$ is the Hopf centre of $A$, i.e., largest central Hopf subalgebra of~$A$.

Recall that a braided fusion subcategory $\cc$ is called {\it nondegenerate} if its M\"uger center is trivial, i.e., $\cc'=\mtr{Vec}$. If $\cc$ is a nondegenerate braided fusion category then by \cite[Corollary~3.11] {DGNO} one has
\begin{gather}\label{cadc}
(\cc_{\rm ad})'=\cc_{\rm pt}.
\end{gather}
Since $\cd''=\cd$ for any fusion subcategory $\cd$ of a nondegenerate braided category $\cc$ we can also write that
\begin{gather}\label{cadc2}
(\cc_{\rm pt})'=\cc_{\rm ad}.
\end{gather}
\begin{Proposition} Let $A$ be a factorizable semisimple Hopf algebra. With the above notations, one has that
\begin{gather*}%\label{chop}
\phi_{R}((A//K(A))^{*})=A'\qquad \text{and}\qquad \phir((A/I)^*)=K(A).
\end{gather*}
\end{Proposition}

\begin{proof} For a quasitriangular Hopf algebra $(A, R)$ it is well known that $\rep(A)$ is nondegenerate if and only if~$A$ is a factorizable Hopf algebra. In this case equation~\eqref{cadc} gives that
 \begin{gather*}
 \rep(A//K(A))'=\rep(A//A').
 \end{gather*}
Then Theorem \ref{main1} gives $\phir((A//K(A))^*)=A'$.

Similarly, equation~\eqref{cadc2} gives that
 \begin{gather*}
 \rep(A//A')'=\rep(A//K(A))
 \end{gather*}
and Theorem~\ref{main1} implies $K(A)=\phi_{R}((A//A')^{*})$ .\end{proof}

\section{Conjugacy classes and M\"uger centralizer}\label{smr}
In this section we will prove Theorem~\ref{mgqtr}. Suppose that $(A,R)$ is a semisimple quasitriangular Hopf algebra. Let as above $V_{0}, V_{1},\dots ,V_{r}$ be a complete set of isomorphism classes of irreducible $A$-modules. Let also $\irr(A)=\{\ch_{0}, \ch_{1}, \dots ,\ch_{r}\}$ be the set of irreducible characters afforded by these modules and $\{E_{0}, \dots ,E_r\}$ be their associated central primitive idempotents of~$A$. Without loss of generality we may suppose that $V_0=\kk$ is the trivial $A$-module and therefore $\ch_0=\eps$ and $E_0=\blam$.

Recall by \cite{EG}, that $\rep(A)$ is a ribbon category with the canonical ribbon element $v=u^{-1}$, where $u:=S\big(R^1\big)R^2$ is the Drinfeld element of $(A,R)$. With respect to the canonical ribbon structure given by this ribbon element, the $S$-matrix of $(A, R)$ has entries
\begin{gather*}
s_{ii'}:=\mtr{tr}_{V_{i}\ot V_{i'}}(Q)=(\ch_{i}\ot \ch_{i'})(Q)=\langle \ch_{i'}, \phi_{R}(\ch_{i})\rangle.
\end{gather*}
It follows from \cite{dgno2} that one has $|s_{ii'}|\leq \ch_{i}(1)\ch_{i'}(1)$ and $V_{i}, V_{i'}$ centralize each other if and only if
$s_{ii'}=\ch_{i}(1)\ch_{i'}(1).$

The Drinfeld map $\phi_{R}\colon C(A)\ra \mtc Z(A) $ is an algebra map and we may suppose as in the introduction that
\begin{gather*}
\phi_{R}(F_{j})=\sum_{i \in \mtc A_{j}}E_{i}
\end{gather*}
for some subset $\ca_j\subseteq \{0, \dots , r\}$.

Without loss of generality we may also suppose that $F_{0}=t$, the idempotent integral of~$A^{*}$. Then $\phir(F_{0})$ is the idempotent integral of $\kda:=\phir(A^*)$ since $\phir(f)\phir(F_{0})= \phir(fF_{0})=f(1)\phir(F_{0})=\eps(\phir(f))\phir(F_{0})$ for any $f \in A^{*}$ and also $\eps(\phir(F_{0}))=F_0(1)=1$.

 Note that the set $\mtc A_{j}$ is empty if and only if $\phi_{R}(F_{j})=0$. Denote by $\mtc J\subseteq \{0,1,\dots, r\}$ the set of all indices $j$ with $\mtc A_{j}$ not a empty set. Since $\phi_{R}(1)=1$ we obtain in this way a partition for the set of indices of all irreducible representations $\{0,1, \dots, r\}=\bigsqcup\limits_{j\in \mtc J}\mtc A_{j}$.

For any index $0\leq i\leq r$ we denoted by $m(i)$ the unique index $j \in \mtc J$ such that $i \in \mtc A_{j}$. Therefore in this way we obtain a unique function
\begin{gather*}
m\colon \ \{0, 1,\dots, r\}\ra \mtc J
\end{gather*}
with the property that $E_i\phir(F_{m(i)})\neq 0$ for all $i \in \{0, 1,\dots, r\}$.

Recall from Section~\ref{lkern} the definition of the left kernel of an $A$-module.
\begin{Lemma}Let $(A, R)$ be a quasitriangular Hopf algebra and $V_{i}$, $V_{i'}$ be two irreducible $A$-representations. Then, with the above notations the following assertions are equivalent:
\begin{enumerate}\itemsep=0pt
\item[$1)$] $V_{i}$ and $V_{i'}$ centralize each other in $\rep(A)$,
\item[$2)$] $\ch_{i'}F_{m(i)}=\ch_{i'}(1)F_{m(i)}$,
\item[$3)$] $\ch_{i}F_{m(i')}=\ch_{i}(1)F_{m(i')}$,
\item[$4)$] $\cc^{m(i)}\subseteq \lker_{A}(V_{i'})$,
\item[$5)$] $\cc^{m(i')}\subseteq \lker_{A}(V_{i})$.
\end{enumerate}
\end{Lemma}

\begin{proof}For any character $\ch \in C(A)$ write as above
$\ch=\sum\limits_{j=0}^{r}\al_{\ch,j}F_{j}$.
Then one has that
\[
\phi_{R}(\ch)=\sum\limits_{j=0}^{r}\al_{\ch,j}\phi_{R}(F_{j})=\sum\limits_{j=0}^{r}\al_{\ch,j}\bigg(\sum\limits_{s \in \mtc A_{j}}E_{s}\bigg).
\]
With these formulae note that
\begin{gather*}%\label{sii}
s_{ii'}=\langle \ch_{i}, \phi_{R}(\ch_{i'})\rangle =\bigg\langle \ch_{i}, \sum_{j\in \mtc J}\sum_{s\in \mtc A_{j}}\al_{\ch_{i'}, j}E_{s}\bigg\rangle =\ch_{i}(1)\al_{\ch_{i'}, m(i)},
\end{gather*}
where $m(i)$ as above, is the unique index $j \in J$ with $i \in \mtc A_{j}$. Therefore we see that~$V_{i}$ centralize~$V_{i'}$ if and only if
\begin{gather}\label{centralize}
\al_{\ch_{i'}, m(i)}=\ch_{i'}(1) \ \iff \ \ch_{i'}F_{m(i)}=\ch_{i'}(1) \ \iff \ \cc^{m(i)}\subseteq \lker_{A}(V_{i'}).
\end{gather}

The equivalence of assertions (2) and (4) follows from \cite[Theorem~3.6]{CW5}. The rest of the equivalences follow from the symmetry property of the centralizer.
\end{proof}

\begin{Remark}The above lemma also shows that if $V_{i}$ centralizes $V_{i'}$ then $V_{i}$ centralize all $V_{i''}$ with $i''\in \mtc A_{m(i')}$.
\end{Remark}
Next theorem is a generalization of Theorem~\ref{mgqtr}.

\subsection{Proof of Theorem \ref{mgqtr}}

\begin{proof}Using the previous lemma we have the following equalities:
\begin{align*}
\co(\rep(A//L)')& =\bigcap_{M\in \irr(A//L)}\co(\langle M\rangle ')
 =\bigcap_{M\in \irr(A//L)}\big\{V_{i}\,|\, \cc^{m(i)}\subseteq \lker_{A}(M)\big\}\\
& =
\bigg\{V_{i}\,|\,\cc^{m(i)}\subseteq \bigcap_{M\in \irr(A//L)} \lker_{A}(M)\bigg\}.
\end{align*}
On the other hand by equation~\eqref{brint} one has that
\begin{gather*}
\cap_{M\in \irr(A//L)} \lker_{A}(M)=L
\end{gather*}
and therefore $\co(\rep(A//L)')=\big\{\ch_{i}\,|\,\cc^{m(i)}\subseteq L\big\}$.
\end{proof}

Another description of the simple objects of $\rep(A//L)'$ is given in the following:
\begin{Proposition}\label{evlamqtr}If $(A, R)$ is a semisimple quasitriangular Hopf algebra then
\begin{gather*}
\co(\rep(A//L)')=\{\ch_{i}\,|\, m(i) \in \mtc I_{L}\}
=\{\ch_{i}\,|\, F_{m(i)}(\blam_{L})\neq 0\}\end{gather*}
for any left normal coideal subalgebra $L$ of $A$.
\end{Proposition}

\begin{proof} By Theorem \ref{mgqtr} one has that $\co(\rep(A//L)')=\big\{\ch_{i}\,|\, \cc^{m(i)}\subseteq L\big\}$, i.e., $\co(\rep(A//L)')=\{\ch_{i}\,|\, m(i)\in \mtc I_L \}$. On the other hand by Lemma~\ref{charval} one has that $\cc^{m(i)}\subseteq L$ if and only if $F_{m(i)}(\blam_L)\neq 0$.
\end{proof}

Corollary~\ref{putgh1} together with Theorem~\ref{mgqtr} gives the following:
\begin{Theorem}\label{puttgh} Let $(A, R)$ be a quasitriangular semisimple Hopf algebra and~$L$ be a left normal coideal subalgebra of $A$. For an irreducible representation~$V_{i}$ of $A$ we have that the following assertions are equivalent:
\begin{enumerate}\itemsep=0pt
\item[$1)$] $V_{i}\in \co(\rep(A//L)')$,
\item[$2)$] $\cc^{m(i)}\subseteq L $,
\item[$3)$] $\phi_{R}((A//L)^{*})\subseteq \lker_{A}(V_{i})$.
\end{enumerate}
\end{Theorem}

\subsection[Description of $\kda$]{Description of $\boldsymbol{\kda}$}

As above denote by $\kda:=\phir(A^*)$ the image of the Drinfeld map.
\begin{Proposition}\label{mcenter} Suppose that $(A,R)$ is a quasitriangular Hopf algebra and $\phi_{R}(F_0) =\sum\limits_{i \in \mtc A_{0}}E_{i}$ where $F_0=t$ is the idempotent integral of $A^{*}$. Then
\begin{enumerate}\itemsep=0pt
\item[$1)$] $\rep(A)'=\rep(A//\kda)$,
\item[$2)$] $\irr(A//\kda)=\{\ch_{i}\,|\,i \in \mtc A_{0}\}$,
\item[$3)$] $\kda=\bigoplus\limits_{j \in \mtc J}{\cc^{j}}$.
\end{enumerate}
\end{Proposition}
\begin{proof}By Theorem \ref{main1} one has that
\begin{gather*}
\rep(A)'=\rep(A//\kk)'=\rep(A//\phir(A^*))=\rep(A//\kda).
\end{gather*}
On the other hand, Theorem~\ref{mgqtr} gives the following equality
\begin{gather*}
\co(\rep(A)')=\co(\rep(A//\kk)')=\big\{\ch_{i}\,|\,\cc^{m(i)}\subseteq \kk\big\}=\big\{\ch_{i}\,|\,\cc^{m(i)}= \kk\big\}.
\end{gather*}
It is easy to see that $\cc^{m(i)}= \kk$ if and only if $m(i)=0$. Indeed, $\cc^{m(i)}=\mtc F(F_{m(i)}A^*)= \blam \lh F_{m(i)}A^*$ and $\kk=\mtc F(F_{0}A^*)$. Since $\mtc F$ is bijective the statement follows. Therefore $\co(\rep(A//\kda))=\co(\rep(A)')=\{\ch_{i}\,|\,m(i)=0\}=\{\ch_{i}\,|\,i \in \mtc A_{0}\}$.

For the last item note that $\co(\rep(A))=\co(\rep(A//\kda)')=\big\{\ch_{i}\,|\,\cc^{m(i)}\subseteq \kda\big\}$. This implies that $\kda=\bigoplus_{j \in \mtc J}{\cc^{j}}$.
\end{proof}

\subsection{Proof of Corollary \ref{main2}}
This is now a particular case of Theorem~\ref{mgqtr}.

\begin{proof} Note that if $A$ is factorizable then $\phir$ is bijective and every set $\mtca_j$ is a~singleton. Moreover $\phir(F_0)=E_0$, the integral of~$A$ in this case. Then, without loss of generality, after a permutation of the indices, we may suppose $\phir(F_i)=E_i$ and therefore $m(i)=i$ for all $0\leq i \leq r$. Then the statement of Theorem~\ref{mgqtr} becomes Theorem~\ref{main2}.
\end{proof}

For the rest of this subsection we suppose that $A$ is a semisimple factorizable Hopf algebra. As explained above without loss of generality we may also assume $\phir(F_i)=E_i$ and therefore that the function $m\colon \{0, 1,\dots, r\}\ra \{0, 1,\dots, r\}$ is the identity map.
 \begin{Proposition} Suppose that $(A,R)$ is a semisimple factorizable Hopf algebra. Then for any irreducible $A$-module $V_{i}$ one has that
\begin{gather*}
\lker_{A}(V_{i})=\bigoplus_{\{i' \,|\, V_{i'}\;\text{centralize}\; V_i\}}{\cc^{i'}}.
\end{gather*}
\end{Proposition}
\begin{proof} Since $\mtc J=\{0, 1,\dots, r\}$ in this case, by equation~\eqref{centralize} and \cite[Theorem 3.6]{CW5} one has that
\begin{gather*}
\ch_{i}F_{i'}=\ch_{i}(1)F_{i'} \ \iff \ {\cc^{i'}}\subseteq \lker_{A}(V_{i}) \ \iff \ V_{i}\;\text{centralizes}\; V_{i'}.
\end{gather*}
It follows that in this case one has
\begin{gather*}
\lker_{A}(V_{i})=\bigoplus_{\cc^{i'}\subseteq\lker_{A}(V_{i})}\cc^{i'}=\bigoplus_{\{i'\,|\, V_{i}\;\text{centralizes}\; V_{i'}\}}\cc^{i'}.\tag*{\qed}
\end{gather*}\renewcommand{\qed}{}
\end{proof}

\begin{Remark}\quad
\begin{enumerate}\itemsep=0pt
\item From the previous proposition, in the case of a semisimple factorizable Hopf algebra one can deduce that for any two irreducible characters $\chi_{i}$ and $\chi_{i'}$ one has
$\mtc C^{i'}\subseteq \lker_{A}(V_{i})\iff\mtc C^{i}\subseteq \lker_{A}(V_{i'})$
\item Recall that in \cite[Theorem 1.4]{mathz} it is shown that
\begin{gather*}
\co(\rep(A//K)')=\{\ch_{i}\,|\, F_{i}(\blam_{K})\neq 0\}
\end{gather*}
for any normal Hopf subalgebra $K$ of a factorizable Hopf algebra~$A$. Note that Proposition~\ref{evlamqtr} generalizes the above result from normal Hopf subalgebras~$K$ to left normal coideal subalgebras $L$ of $A$. It also drops the factorizability assumption on~$A$.
\end{enumerate}
\end{Remark}

Define $C_{V_{i}}:=C_{\ch_{i}}\subset A^{*}$ as the subcoalgebra of $A^{*}$ generated by $\ch_{i}$. By \cite[Lemma~4.2(i)]{CW5}, in the factorizable case one has that $\phi_{R}(C_{V_{i}})=\mtc C_{i}.$ for all $0\leq i \leq r$.
\begin{Remark} Let $A$ be a semisimple factorizable Hopf algebra and~$L$ a left normal coideal subalgebra of~$A$. If $\cd:=\rep(A//L)$ then
$(A//L)^{*}=\bigoplus\limits_{\ch_{j} \in \co(\cd)}C_{V_{j}}$
and $L^\circ=\phi_{R}((A//L)^{*})=\bigoplus\limits_{\ch_{j} \in \co(\cd)}\phi_{R}(C_{V_{j}})$, i.e., $L^\circ=\bigoplus\limits_{\ch_{j} \in \co(\cd)}\cc^{j}$. This gives another proof for Theorem~\ref{main2} in the case of a factorizable Hopf algebra since in this case $\cd''=\cd$.
\end{Remark}

\section[Example $H_8$]{Example $\boldsymbol{H_8}$}\label{h8}

In this section we compute the centralizer of any fusion subcategory of the quasi-triangular Hopf algebra $H_8$, the unique semisimple non-trivial Hopf algebra of dimension~$8$. We note that the category of representations~$\rep(H_8)$ is a braided Tambara--Yamagami category and therefore~$\rep(H_8)\subset \rep(D(H_8))$. The $S$-matrix of the center of a Tambara--Yamagami was computed in~\cite{gnn}. Using this one can describe completely the centralizer of any fusion subcategory of~$\rep(H_8)$. However, we decided to include this example here to illustrate how Theorem~\ref{mgqtr} can be applied in a concrete example.

The eight-dimensional semisimple Hopf algebra (see \cite{kp-66, ma-6-8}) is generated by $\{x,y,z\}$ subject to the relations
\begin{gather*}
x^2=y^2=z^2=1,\qquad xz=zx,\qquad zy=yz, \qquad xyz=yx.
\end{gather*}
The comultiplication is given by
\begin{gather}\label{d1}\Delta(x)=xe_0\ot x+xe_1\ot y, \qquad \Delta(y)=ye_1\ot x+ye_0\ot y,\qquad \Delta(z)=z\ot z,
\end{gather}
and the counit by $\eps(x)=\eps(y)=\eps(z)=1.$ The antipode has the formulae
\begin{gather*}%\label{ant}
S(x)=xe_0+ye_1,\qquad S(y)=xe_1+ye_0,\qquad S(z)=z.
\end{gather*}
Based on equation~\eqref{d1} one can compute that
\begin{gather*}%\label{d2}
\Delta(xy)=xye_0\ot xy+xye_1\ot yx,\qquad \Delta(yx)=yxe_0\ot yx+yxe_1\ot xy.
\end{gather*}
It can also be checked that
\begin{gather*}%\label{d3}
\Delta(xz)=xe_0\ot xz-xe_1\ot yz,\qquad \Delta(yz)=ye_0\ot yz-ye_1\ot xz.
\end{gather*}
Since $z$ is a central element in $H_8$, there are two central orthogonal idempotents:
\begin{gather*}
e_0=\frac{1}{2}(1+z),\qquad e_1=\frac{1}{2}(1-z).
\end{gather*}
\subsection[Dual basis of $H^*$]{Dual basis of $\boldsymbol{H^*}$}

$H_8$ has a $\kk$-linear basis given by the set of elements $\{1, x, y, z, {xy}, {yx}, {xz}, {yz}\}$. We consider its linear dual basis on $H_8^*$ given by $\{p_1, p_x, p_y, p_z, p_{xy}, p_{yx}, p_{xz},p_{yz}\}$.

It is easy to see that the idempotent integrals of $H_8$ and $H_8^*$ are given by
\begin{gather*}%\label{intgerals}
\blam=\frac{e_0}{4}(1+x+y+xy)=\frac{1}{8}(1+x+y+xy+z+zx+zy+yx),\qquad \lam=p_1.
\end{gather*}

\subsection[$H_8^*$ representations]{$\boldsymbol{H_8^*}$ representations}

Since $H_8$ is a self dual Hopf algebra~\cite{alaoui} it has also four-1-dimensional representations given by the group like elements of $H_8$ and a~2-dimensional representation.
One has that
\begin{gather*}
G(H_8)=\{1, g_1, g_2, z\},
\end{gather*}
with $g_1=xy(e_0+ie_1)$, $g_2=xy(e_0-ie_1)$. It can be easily checked that $g_1g_2=z$, $zg_i=g_i$ and $g_i^2=1$. Moreover, the set of central grouplike elements of~$H_8$ is given by $\bar G(H_8)=\{1, z\}$.

\subsubsection{On the 2-dimensional comodule}

From equation~\eqref{d1} one can compute that
\begin{alignat*}{3}
& \Delta(xe_0)=xe_0\ot xe_0+xe_1\ot ye_1,\qquad && \Delta(xe_1)=xe_0\ot xe_1+xe_1\ot ye_0,&
\\
& \Delta(ye_0)=ye_0\ot ye_0+ye_1\ot xe_1,\qquad && \Delta(ye_1)=ye_0\ot ye_1+ye_1\ot xe_0.&
\end{alignat*}
This shows that $W=\kk w_1\oplus \kk w_2$ is a left $H_8$-comodule with
the comodule structure given by
\begin{gather*}
\rho(w_1)=xe_0\ot w_1+xe_1\ot w_2,\qquad \rho(w_2)=ye_1\ot w_1+ye_0\ot w_2.
\end{gather*}
Thus
\begin{gather}\label{comsd}
H_8=\kk 1\oplus \kk g_1 \oplus \kk g_2\oplus \kk z \oplus(\kk xe_0\oplus \kk ye_1 )\oplus (\kk ye_0 \oplus \kk xe_1)
\end{gather}
is a decomposition of $H_8$ into simple left $H_8$-comodules. We denote $M_0=\kk 1$, $M_1=\kk g_1$, $ M_2=\kk g_2$, $M_3=\kk z$ the four one-dimensional $H_8^*$-modules. Moreover, $M_4:= (\kk xe_0\oplus \kk ye_1 )$ is an irreducible $H_8^*$-module and $M_4\simeq M_5:=(\kk ye_0 \oplus \kk xe_1)$.

\subsection{Fourier transform} Based on the comultiplication formulae one can compute the Fourier transform
\begin{gather*}
\mtc F\colon \ H_8^*\ra H_8,\qquad f \mapsto \dimk (H_8)f\rh \blam.
\end{gather*}
After some computations it follows that under $\mtc F$ one has
\begin{gather*}
p_1 \mapsto 1, \qquad p_z\mapsto z,\;p_x\mapsto xe_0+ye_1, \qquad p_{xz}\mapsto xe_0-ye_1
\end{gather*}
and
\begin{gather*}
p_y\mapsto ye_0+xe_1, \qquad p_{yz}\mapsto ye_0-xe_1,\qquad p_{xy}\mapsto yx, \qquad p_{yx}\mapsto xy.
\end{gather*}

\subsection[Irreducible $H_8$-modules and their characters]{Irreducible $\boldsymbol{H_8}$-modules and their characters}

$H_8$ has four-1-dimensional modules $V_0$, $V_1$, $V_2$, $V_3$ and a 2-dimensional irreducible module~$V$. The action of the generators on these modules is given as follows.

For $V_0$, the trivial $H$-module the action is given by $xv=yv=zv=v$ and therefore the character is given by
\begin{gather*}
\ch_0=p_1+p_z+p_x+p_y+p_{xy}+p_{yx}+p_{xz}+p_{yz}.
\end{gather*}
For $V_1$ the action is given by $xv=-v$, $yv=v$, $zv=v.$ The character $\ch_1$ of $V_1$ has $\ch_1(1)=\ch_1(z)=\ch_1(y)=\ch_1(yz)=1$, $\ch_1(x)=\ch_1(yx)=\ch_1(xy)=\ch_1(xz)=-1$. Thus
\begin{gather*}
\ch_1=p_1+p_z-p_x+p_y-p_{xy}-p_{yx}-p_{xz}+p_{yz}.
\end{gather*}
For $V_2$ the action is given by $xv=-v$, $yv=-v$, $zv=v$ and the character is given by
\begin{gather*}
\ch_2=p_1+p_z-p_x-p_y+p_{xy}+p_{yx}-p_{xz}-p_{yz}.
\end{gather*}
 For $V_3$ one has $xv=v$, $yv=-v$, $zv=v$ and the character is given by
\begin{gather*}
\ch_3=p_1+p_z+p_x-p_y-p_{xy}-p_{yx}+p_{xz}-p_{yz}.
\end{gather*}
For $V_4=V$, the 2-dimensional simple module, if $V=\kk v_1\oplus \kk v_2$ then the action of generators given by
\begin{gather*}
xv_1=v_1,\qquad yv_1=v_2,\qquad zv_1=-v_1,\qquad xv_2=-v_2,\qquad yv_2=v_1,\qquad zv_2=-v_2.
\end{gather*}
The character $\ch_4$ of $V$ has the values $\ch_4(1)=2 \ch_4(z)=-2$, $\ch_4(x)=\ch_4(y)= \ch_4(yx)=\ch_4(xy)=\ch_4(xz)=\ch_4(zx)=0$. This gives that
\begin{gather*}
\ch_4=2p_1-2p_z.
\end{gather*}
\subsubsection{Multiplication of the characters}It can be easily checked that $\ch_4^2=\sum\limits_{i=0}^3\ch_i$, and $\ch_i\ch_4=\ch_4\ch_i=\ch_i$ for all $0\leq i \leq 3$. Moreover $\ch_i\ch_j=\ch_k$ if $\{i,j,k\}=\{1,2,3\}$.

\subsection{Central primitive idempotent of the character ring}
Based on the above multiplication, the central primitive idempotents of $C(H_8)$ can be computed as follows
\begin{gather*} F_0=\frac{1}{8}(\eps+\ch_1+\ch_2+\ch_3+2\chi_4)=p_1 \qquad \text{is the integral of $H_8^*$},\\
F_1=\frac{1}{4}(\eps+\ch_1-\ch_2-\ch_3)=p_y+p_{yz}, \qquad F_2=\frac{1}{4}(\eps-\ch_1+\ch_2-\ch_3)=p_{xy}+p_{yx},
\\ F_3=\frac{1}{4}(\eps-\ch_1-\ch_2+\ch_3)=p_x+p_{xz},\qquad F_4=\frac{1}{8}(\eps+\ch_1+\ch_2+\ch_3-2\chi_4)=p_z.
\end{gather*}
Note that $F_4$ is the central primitive idempotent of $H_8^*$ attached to the central grouplike element $z\in \bar G(H_8):=G(H_8)\cap \mtc Z(H_8)$.

\subsection[Conjugacy class sums of $H_8$]{Conjugacy class sums of $\boldsymbol{H_8}$} Using the above formulae for the central primitive idempotents since $ {C}_i=8\big(\blam\lh F_i\big)$ it follows that
\begin{gather*}
 {C}_0=8(\blam\lh F_0)=1,\qquad {C}_1=8(\blam\lh F_1)=2ye_0,\\
 {C}_2=8(\blam\lh F_2)=xy+yx=g_1+g_2,
\end{gather*}
and
\begin{gather*}
 {C}_3=8(\blam\lh F_3)=2xe_0,\qquad {C}_4=8(\blam\lh F_4)=z.
\end{gather*}

\subsection{Description of the adjoint action}
Using the antipode formulae one can see that the adjoint action of $H_8$ on itself can be given by
\begin{gather*}%\label{adjc}
x.a=xax,\qquad y.a=yay,\qquad z.a=a,\qquad \text{for all} \ a\in H_8.
\end{gather*}
\subsection[Conjugacy classes of $H_8$]{Conjugacy classes of $\boldsymbol{H_8}$}

Recall that the conjugacy classes are simple $D(H_8)$-modules~\cite{zind}. We rewrite the decomposition of $H_8$ into left $H_8$-comodules from equation~\eqref{comsd} as follows
\begin{gather*}
H_8=\kk 1\oplus \kk z\oplus (\kk g_1\oplus \kk g_2)\oplus (\kk xe_0\oplus \kk ye_1) \oplus (\kk xe_1\oplus \kk ye_0).
\end{gather*}
Moreover, by the above description of the left adjoint action it can be checked each of the above five subspaces is closed under the adjoint action of~$H_8$.

Clearly, the first two subspaces are irreducible $D(H_8)$-modules being one-dimensional. Since $ C_0\in \kk 1$ and $ C_4=z\in \kk z$ we deduce that
\begin{gather*}
\cc^0=\kk 1,\qquad \cc^4=\kk z.
\end{gather*}
It can be checked directly that the third $D(H_8)$-module is an irreducible $D(H_8)$-module since
\begin{gather*}
x.g_1=g_2,\qquad x.g_2=g_1,\qquad y.g_1=g_2,\qquad y.g_2=g_1.
\end{gather*}
Since $ C_2=g_1+g_2=xy+yx$ it follows that
\begin{gather*}
\cc^2=\kk g_1\oplus \kk g_2.
\end{gather*}
 Similarly one can check that the simple $H_8$-comodule $M_4=\kk xe_0\oplus \kk ye_1$ is an irreducible $D(H_8)$-module and since $ C_3=2xe_0=x+xz\in M_4$ we can say that
\begin{gather*}
\cc^3=\kk xe_0\oplus \kk ye_1.
\end{gather*}
By the same argument
\begin{gather*}
\cc^1=\kk xe_1\oplus \kk ye_0.
\end{gather*}

\subsection[Presentation of the central primitive idempotents of $H_8$]{Presentation of the central primitive idempotents of $\boldsymbol{H_8}$}

The associated central primitive idempotents of $V_0$, $V_1$, $V_2$, $V_3$, $V$ can be computed as
\begin{gather*}
E_0=\frac{e_0}{4}(1+x+y+xy)=\blam
\end{gather*}
is the central primitive idempotent of $\ch_0$, i.e., the idempotent integral~$\blam_{H_8}$ of~$H_8$,
\begin{gather*}
E_1=\frac{e_0}{4}(1-x+y-xy)
\end{gather*}
is the central primitive idempotent of~$\ch_1$,
\begin{gather*}
E_2=\frac{e_0}{4}(1-x-y+xy)
\end{gather*}
is the central primitive idempotent of $\ch_2$,
\begin{gather*} E_3=\frac{e_0}{4}(1+x-y-xy)
\end{gather*}
is the central primitive idempotent of $\ch_3$,
\begin{gather*}
E_4=e_1
\end{gather*}
is the central primitive idempotent of $\ch_4$.

\subsection[The $R$-matrix and the Drinfeld map]{The $\boldsymbol{R}$-matrix and the Drinfeld map}

It is well-known that $H_8$ is a semisimple quasitriangular Hopf algebra~\cite{ma-6-8} with the $R$-matrix given by
\begin{gather*}
R=\frac{1}{2}(1\ot 1+g_2\ot 1+ 1\ot g_1-g_2\ot g_1).
\end{gather*}
It follows that
\begin{gather*}
Q=R_{21}R = \frac{1}{4}(1\ot (1+z+g_1+g_2) + g_1\ot (1+g_1-g_2-z)\\
\hphantom{Q=R_{21}R =}{}+ g_2\ot (1-g_1+g_2-z)+z\ot (1-g_1-g_2+z).
\end{gather*}

\subsubsection[The Drinfeld map on $\kk$-linear basis]{The Drinfeld map on $\boldsymbol{\kk}$-linear basis}
One can compute that the Drinfeld map $\phir\colon H_8^*\ra H_8$, $f\mapsto f\big(Q^1\big)Q^2$ is given by
\begin{gather*}
p_1\mapsto \Lam_{\kk G}=\frac{1}{4}(1+g_1)(1+g_2)=e_{0,0},\qquad p_z\mapsto \frac{1}{4}(1-g_1)(1-g_2)= e_{11},
\\
p_{xy}\mapsto \frac{1}{4}(1+ig_1-ig_2-z),\qquad p_{yx}\mapsto \frac{1}{4}(1-ig_1+ig_2-z),\qquad p_x, p_y, p_{xz}, p_{yz}\mapsto 0.
\end{gather*}

\subsubsection[The Drinfeld map given on the central primitive idempotents of the character ring]{The Drinfeld map given on the central primitive idempotents\\ of the character ring}

The following formulae hold for the Drinfeld map:
\begin{gather*}
\phir (F_0)=E_0+E_2,\qquad \phir (F_1)=\phir (F_3)=0,\qquad \phir (F_2)=E_4,\qquad \phir (F_4)=E_1+E_3.
\end{gather*}

\subsubsection[The map $m$ and the subset $\mtc J$]{The map $\boldsymbol{m}$ and the subset $\boldsymbol{\mtc J}$}

It follows that in this case the subset
$\mtc J=\{0,2,4\}$ and the function $m\colon \{0,1,2,3,4\}\ra \mtc J$ is given by
\begin{gather} \label{m}
m(0)=m(2)=0, \qquad m(1)=m(3)=4, \qquad m(4)=2.
\end{gather}
Therefore $\phir(H^*_8)=\kk G(H_8)$ and by Proposition \ref{mcenter} one has
\begin{gather*}
\co(\rep(H_8)')=\{\ch_i\,|\, m(i)=0\}=\{\ch_0, \ch_2\}.
\end{gather*}
\subsection[Description of the centralizer based on the function $m$]{Description of the centralizer based on the function $\boldsymbol{m}$}

In this subsection we describe the centralizer of any fusion subcategory of $\rep(H_8)$. For the entire category, $\rep(H_8)$ this was done above.

There are two normal coideal subalgebras of $H_8$ which are not Hopf subalgebras,
\begin{gather*}
L_1=\kk\langle 1, z, ye_0, xe_1\rangle =\lker_H(V_1)
\end{gather*}
and
\begin{gather*}
L_3=\kk\langle 1, z, xe_0, ye_1\rangle =\lker_H(V_3).
\end{gather*}
One has
\begin{gather*}
L_1=\cc^0\oplus \cc^1\oplus \cc^4,\qquad L_2=\cc^0\oplus \cc^3\oplus \cc^4.
\end{gather*}
Theorem \ref{charofim} gives that $\rep(H_8//\lker_{H_8}(V_i))\simeq \langle V_i\rangle$, the fusion subcategory of $\rep(H_8)$ generated by $V_i$. Using Theorem \ref{mgqtr} it follows
\begin{gather*}
\rep(A//L_1)'=\big\{\ch_i\,|\,\cc^{m(i)}\subseteq L_1\big\}=\{\ch_i\,|\,m(i)\in \{0,1,4\}\}.
\end{gather*}
By description of the function $m$ from equation~\eqref{m} it follows that
\begin{gather*}
\co(\langle V_1\rangle')=\{\ch_0, \ch_1, \ch_2, \ch_3\}.
\end{gather*}
Similarly, one has that
\begin{gather*}
\co(\langle V_3\rangle')=\rep(A//L_3)'=\big\{\ch_i\,|\,\cc^{m(i)}\subseteq L_3\big\}=\{\ch_i\,|\, m(i)\in \{0,1,3\}\}.
\end{gather*}
By description of the function $m$ from equation \eqref{m} it follows that
\begin{gather*}
\co(\langle V_3\rangle ')=\{\ch_0, \ch_1, \ch_2, \ch_3\}.
\end{gather*}
Note that for the central linear character $\ch_2\in \mtc Z(H_8^*)$ one has
\begin{gather*}
L_2:=\lker_H(V_2)=\kk G(H_8)=\cc^0\oplus \cc^2\oplus \cc^4.
\end{gather*}
As above, using Theorem~\ref{mgqtr} it follows
\begin{gather*}
\co(\langle V_2\rangle')=\{\ch_i\,|\,m(i)\in \{0,2,4\}\}=\co(\rep(H_8)).
\end{gather*}
On the other hand since $\langle V_4\rangle =\rep(H_8)$ Theorem~\ref{charofim} gives that
$L_4=\lker_{H_8}(V)=\kk=\cc^0.$ Thus, using again Theorem~\ref{mgqtr} it follows that
\begin{gather*}
\co(\langle V_4\rangle')=\{\ch_i\,|\,m(i)=0\}=\{\ch_0, \ch_2\}=\co(\rep(H_8)').
\end{gather*}
All the four one-dimensional centralize each other. $V$ centralizes only $\ch_2$.

\subsubsection{On the first commutator and adjoint subcategory}
One has that the first commutator of $H_8$ is given by $H_8'=\kk \langle 1, z\rangle$ and moreover, for this Hopf algebra $\ccad=\ccpt$. Thus in this case
\begin{gather*}
\ccad'=\langle V_1, V_2, V_3\rangle'=\bigcap_{i=0}^3 \langle V_i\rangle'=\ccpt.
\end{gather*}

\subsection*{Acknowledgements}

This research was supported by a grant of Romanian Ministry of Research and Innovation CNCS-UEFISCDI, Project No. PN-III-P4-ID-PCE-2016-0157, within PNCDI III.

\pdfbookmark[1]{References}{ref}
\LastPageEnding

\end{document}